\magnification=1200
\input amssym.def
\input amssym

\font\sectbf=cmbx10 scaled\magstep1
\newskip\ttglue
\font\eightrm=cmr8
\font\eightbf=cmbx8
\font\eightit=cmti8
\font\eighti=cmmi8
\font\eightsl=cmsl8
\font\eightsy=cmsy8
\font\eighttt=cmtt8
\font\sixrm=cmr6
\font\sixi=cmmi6
\font\sixsy=cmsy6

\def\eightpoint{\def\rm{\fam0\eightrm}
  \textfont0=\eightrm \scriptfont0=\sixrm \scriptscriptfont0=\fiverm
  \textfont1=\eighti  \scriptfont1=\sixi  \scriptscriptfont1=\fivei
  \textfont2=\eightsy \scriptfont2=\sixsy \scriptscriptfont2=\fivesy
  \textfont3=\tenex   \scriptfont3=\tenex \scriptscriptfont3=\tenex
  \textfont\itfam=\eightit  \def\it{\fam\itfam\eightit}%
  \textfont\slfam=\eightsl  \def\sl{\fam\slfam\eightsl}%
  \textfont\ttfam=\eighttt  \def\tt{\fam\ttfam\eighttt}%
  \textfont\bffam=\eightbf  \def\bf{\fam\bffam\eightbf}%
   \scriptscriptfont\bffam=\fivebf  \def\bf{\fam\bffam\eithgtbf}%
  \tt \ttglue=.5em plus.25em minus.15em
  \normalbaselineskip=9pt
  \setbox\strutbox=\hbox{\vrule height7pt depth2pt width0pt}%
  \let\sc=\sixrm  \let\big=\eightbig  \normalbaselines\rm}

\def\minusdot{\hbox{\hskip5pt\raise3pt\hbox{$\cdot$}\hskip-5pt$-$}\hskip2pt}
\font\tenrsfs=rsfs10
\font\sevenrsfs=rsfs7
\font\fiversfs=rsfs5
\newfam\rsfsfam
\textfont\rsfsfam=\tenrsfs
\scriptfont\rsfsfam=\sevenrsfs
\scriptscriptfont\rsfsfam\fiversfs

\font\Huge=cmr10 scaled\magstep4

\catcode`\@=11
\newdimen\ex@
\ex@.2326ex
\def\setboxz@h{\setbox\z@\hbox}
\def\leftrightarrowfill@#1{\setboxz@h{$#1-\m@th$}\ht\z@\z@
  $#1\m@th\mathord\leftarrow\mkern-6mu\cleaders
  \hbox{$#1\mkern-2mu\box\z@\mkern-2mu$}\hfill
  \mkern-6mu\mathord\rightarrow$}
\def\overleftrightarrow{\mathpalette\overleftrightarrow@}
\def\overleftrightarrow@#1#2{\vbox{\ialign{##\crcr\leftrightarrowfill@#1\crcr
 \noalign{\kern-\ex@\nointerlineskip}$\m@th\hfil#1#2\hfil$\crcr}}}

\catcode`\@=12
\def\cf{{\rm cf}}

\def\scrp#1{{\fam\rsfsfam\relax#1}}

\def\[{$$}
\def\]{$$}

\def\qedbox{\vbox{\hrule\hbox{\vrule\kern3pt
\vbox{\kern3pt\kern3pt}\kern3pt\vrule}\hrule}}
\def\iff{\hbox{if\hskip .01in f}\ }
\def\underl#1 {\leavevmode\underli #1 }
\def\underli#1 {\ifx&#1\let\next=\relax\unskip
                \else\let\next=\underli \ulinebox{#1} \fi\next}
\def\ulinebox#1{\vtop{\hbox{\strut#1}\hrule}}
\def\today{\ifcase\month\or January\or February\or March\or April\or
May\or June\or July\or August\or September\or October\or November\or
December\fi \space\number\day, \number\year} 

{\obeyspaces\global\let =\ } %

\def\boxtext#1 {%
\lower 3.8pt \hbox{%
\vbox{%
\hrule
\hbox{\strut \vrule$\,$#1$\,$\vrule}%
\hrule
}
}
}
\def\blacktext#1 {%
\lower 3.8pt \hbox{%
\vbox{%
\hrule depth 2pt
\hbox{\strut \vrule width 2pt$\,$#1$\,$\vrule width 2pt}%
\hrule depth 2pt
}
}
}
\def\boxit#1{\vbox{\hrule\hbox{\vrule\kern3pt
    \vbox{\kern3pt#1\kern3pt}\kern3pt\vrule}\hrule}}
\def\bigfirstletter#1#2{{\noindent
    \setbox0\hbox{\Huge #1}\setbox1\hbox{#2}\setbox2\hbox{(}%
    \count0=\ht0\advance\count0 by\dp0\count1\baselineskip
    \advance\count0 by-\ht1\advance\count0 by\ht2
    \dimen1=.5ex\advance\count0 by\dimen1\divide\count0 by\count1
    \advance\count0 by1\dimen0\wd0
    \advance\dimen0 by.25em\dimen1=\ht0\advance\dimen1 by-\ht1
    \global\hangindent\dimen0\global\hangafter-\count0
    \hskip-\dimen0\setbox0\hbox to\dimen0{\raise-\dimen1\box0\hss}%
    \dp0=0in\ht0=0in\box0}#2}
\def\move(#1,#2){\hskip #1mm\raise#2mm}
\def\markem(#1){\hbox{\special{em:point #1}}}

\def\marginnote#1{\par\noindent\hskip-.5in\hbox to.72in{#1\hfil}}

\def\cell{\hbox{\rm c}}
\def\depth{\hbox{\rm Depth}}
\def\d{\hbox{\rm d}}
\def\length{\hbox{\rm Length}}
\def\irr{\hbox{\rm Irr}}

\def\ind{\hbox{\rm Ind}}
\def\t{\hbox{\rm t}}
\def\s{\hbox{\rm s}}
\def\hl{\hbox{\rm hL}}
\def\inc{\hbox{\rm Inc}}
\def\hcof{\hbox{\rm h-cof}}
\def\ult{\hbox{\rm Ult}}
\def\aut{\hbox{\rm Aut}}

\def\id{\hbox{\rm Id}}
\def\sub{\hbox{\rm Sub}}
\def\hd{\hbox{\rm hd}}
\def\intalg{\hbox{\rm Intalg}\,}

\tolerance=2000

\centerline{\sectbf Concerning problems about}
\centerline{\sectbf cardinal invariants on Boolean algebras}

\medskip
\centerline{J. Donald Monk}
\centerline{\today}

\bigskip
\par\noindent
The purpose of these notes is to describe the progress made on the 97 open
problems formulated in the book {\bf Cardinal invariants on Boolean
algebras}, hereafter denoted by [CI]. Although we assume acquaintance with
that book, we give some background for many problems, and state without
proof most results relevant to the problems, as far as the author
knows. Many of the problems have been solved, at least partially, by
Saharon Shelah. For the unpublished papers of Shelah
mentioned in the references, see his archive,
which can currently be obtained at the URL

\medskip
\par\noindent
http://shelah.logic.at/

\medskip
\par\noindent
We follow the notation of Koppelberg [89] and of [CI]. All Boolean algebras
are assumed to be infinite, unless otherwise indicated.

At this time, the problems fully solved are: 2, 4, 10, 12, 16, 18, 20, 22, 24,
34, 35, 36, 41, 42, 46, 47, 48, 51, 52, 55, 56, 60, 62, 76, 80, 92;
those partially solved
are: 1, 8, 11, 13, 14, 21, 25, 26, 27, 32, 33, 37, 45, 50, 54, 73, 77,
78, 79, 81, 93;
still open in the form given in [CI] are: 3, 5, 6, 7, 9, 15, 17, 19,
23, 28, 29, 30, 31, 38, 39, 40, 43, 44, 49, 53, 57, 58, 59, 61, 63,
64, 65, 66, 67, 68, 69, 70, 71, 72, 74, 75, 82, 83, 84--91, 94--97.

Among those which are still unsolved, or partially unsolved, progress
on the following may be possible without too much difficulty: 5, 7, 8,
9, 11, 13, 14, 21, 23, 29, 30, 31, 32, 33, 37, 39, 43, 44, 45,
50, 54, 64, 68, 69, 70, 71, 72, 73, 74, 77, 78, 81, 84, 86, 90, 93, 94--97.

\medskip
\centerline{\bf Problem 1: cellularity of free products}

\medskip
\par\noindent
In its most general form, this problem is to fully describe what happens
to cellularity in passing from BAs $A,B$ to their free product $A\oplus
B$; more precisely, for which cardinals $\kappa$ and $\lambda$ are
there BAs $A,B$ such that $\cell'(A)=\kappa$, $\cell'(B)=\lambda$, and
$\cell'(A\oplus
B)>\max(\kappa,\lambda)$? Here $\cell'(A)$ is the least cardinal
$\mu$ such that $A$ does not have a disjoint subset of size $\mu$.
Recall from the Erd\"os, Tarski theorem (see Koppelberg [89]) that
$\cell'(A)$ is always a regular cardinal.

In case
$\kappa=\lambda$, the problem is equivalent to asking for a BA $A$ such
that $\cell'(A)=\kappa$ and $\cell'(A\oplus A)>\kappa$. This is because of
the elementary fact that
\[(A\times B)\oplus(A\times B)\cong(A\oplus A)\times(A\oplus
B)\times(B\oplus A)\times(B\oplus B).\]
The problem has been worked on extensively by Shelah and Todor\v cevi\'c.
We mention only results which are not superceded by later work.

\medskip
\par\noindent
(1) $\cell'(A\oplus B)\leq(2^{{\rm c}(A)\cdot{\rm c}(B)})^+$.
Thus under GCH, $\cell'(A\oplus B)\leq\max((\cell(A))^{++},(\cell(B))^{++})$.
See [CI].

\medskip
\par\noindent
(2) If $\cell'(A)\leq(2^\kappa)^+$ and $\cell'(B)\leq\kappa^+$, then
$\cell'(A\oplus B)\leq(2^\kappa)^+$.
Under GCH, if $\cell'(A)\leq\kappa^{++}$ and $\cell'(B)\leq\kappa^+$, then
$\cell'(A\oplus B)\leq\kappa^{++}$.
See [CI].

\medskip
\par\noindent
(3) If $\kappa$ is strong limit, $\cell'(A)\leq\kappa^+$, and
$\cell'(B)\leq\cf(\kappa)$, then $\cell'(A\oplus B)\leq\kappa^+$. See [CI].

\medskip
\par\noindent
(4) It is consistent with ZFC that there is a BA $A$ with
$\cell(A)=\omega$ and $\cell(A\oplus A)>\omega$. This is a classical
result of Kurepa. See [CI].

\medskip
\par\noindent
(5) $\hbox{MA}+\neg\hbox{CH}$ implies that if $\cell(A)=\omega=\cell(B)$,
then $\cell(A\oplus B)=\omega$. See [CI].

\medskip
\par\noindent
Note that (4) and (5) take care of the case $\kappa=\lambda=\omega$
of our problem.

\medskip
\par\noindent
(6) For $\kappa$ regular and $>$ $\omega_1$, there is a BA $A$ such that
$\cell'(A)=\kappa^+$ while $\cell'(A\oplus A)>\kappa^+$. See Shelah [91].

\medskip
\par\noindent
(7) There is a BA $A$ such that $\cell'(A)=\aleph_2$ and $\cell(A\oplus
A)>\aleph_2$. See Shelah [97].

\medskip
\par\noindent
(8) If $\kappa$ is singular, then there is a BA $A$ such that
$\cell'(A)=\kappa^+$ while $\cell'(A\oplus A)>\kappa^+$. See Shelah [94].

\medskip
\par\noindent
Note that (6)--(8) take care of our problem when $\kappa=\lambda$ is a
successor cardinal.

\medskip
\par\noindent
(9) If $\kappa$ is weakly compact and $\cell'(A)=\kappa$, then
$\cell'(A\oplus A)=\kappa$. This is obvious from the equivalent
$\kappa\rightarrow(\kappa)^2$ of weak compactness.

\medskip
\par\noindent
(10) If $\kappa$ is inaccessible and has a stationary subset which does
not reflect in any inaccessible, and if $\cf\delta$ is uncountable for
each $\delta\in S$, then there is a BA $A$ with $\cell'(A)=\kappa$ while
$\cell'(A\oplus A)>\kappa$.

\medskip
\par\noindent
This follows from Shelah [94$'$], 4.8, using
Shelah [94$''$]. Recall that a stationary set on $\kappa$ reflects in a
limit ordinal $\alpha<\kappa$ \iff $S\cap\alpha$ is stationary in
$\alpha$.

These two last results come close to taking care of the case
$\kappa=\lambda$ regular limit.

\medskip
\par\noindent
(11) It is consistent that for all BAs $A,B$, if
$\cell'(A)=\aleph_1$ and $\cell'(B)\leq2^{\aleph_0}$, then
$\cell'(A\oplus B)\leq2^{\aleph_0}$. See Shelah [96].

\medskip
\par\noindent
(12) If $\mu$ is a singular strong limit cardinal, $\theta=(2^{{\rm
cf}(\mu)})^+$, and $2^\mu=\mu^+$, then there are BAs $A,B$ such that
$\cell(A)=\mu$, $\cell(B)<\theta$, and $\cell(A\oplus B)=\mu^+$.
Note that the hypothesis holds for every singular
cardinal under GCH.
See
Shelah [00].

\medskip
\par\noindent
This finishes our survey of results concerning Problem 1. Note that the
following are still, evidently, open.

(a) If $\kappa$ is inaccessible but not weakly compact,
and every stationary subset of $\kappa$
reflects in an inaccessible, is there a BA $A$ such that
$\cell'(A)=\kappa$ while $\cell'(A\oplus A)>\kappa$?

(b) Describe fully the situation of the general problem when
$\kappa\not=\lambda$.

\bigskip
\centerline{\bf Problem 2: cellularity and ultraproducts}

\bigskip
\par\noindent
This problem was solved by Magidor, Shelah [98]: It is consistent that
there is an infinite set $I$, a system $\langle A_i:i\in I\rangle$ of
infinite BAs, and an ultrafilter $F$ on $I$ such that $\cell(\prod_{i\in
I}A_i/F)<|\prod_{i\in I}\cell(A_i)/F|$.

\vfill
\eject
\centerline{\bf Problem 3: the subalgebra-cellularity relation problem}

\medskip
\par\noindent
This general question is still open; no purely cardinal number
characterization of $\cell_{\rm Sr}$ is known.

\bigskip
\centerline{\bf Problem 4: a particular subalgebra-cellularity relation
problem}

\medskip
\par\noindent
This problem was solved negatively by the following theorem in Monk,
Nyikos [97]:

\medskip
\par\noindent
{\it For every infinite cardinal $\kappa$ and every BA $A$, if
$\cell(A)\geq\kappa^{++}$ and $(\kappa,\kappa^{++})\in\cell_{\rm Sr}(A)$,
then $(\kappa^+,\kappa^{++})\in\cell_{\rm Sr}(A)$.}

\medskip
\par\noindent
So in particular, if $\cell(A)=\omega_2$ and
$(\omega,\omega_2)\in\cell_{\rm Sr}(A)$, then
$(\omega_1,\omega_2)\in\cell_{\rm Sr}(A)$, which solves problem 4.

\bigskip
\centerline{\bf Problem 5: two more particular subalgebra-cellularity
relation problems}

\medskip
\par\noindent
These two problems are still open.

\bigskip
\centerline{\bf Problem 6: the homomorphism-cellularity relation problem}

\medskip
\par\noindent
This general problem is still open.

\bigskip
\centerline{\bf Problem 7: two particular homomorphism-cellularity
relation problems}

\medskip
\par\noindent
These two problems are still open.

\bigskip
\centerline{\bf Problem 8: four particular homomorphism-cellularity
relation problems}

\medskip
\par\noindent
Two of these have been solved. Theorem 4 in Monk, Nyikos [97] says:

\medskip
\par\noindent
{\it Suppose that $\omega\leq\rho\leq\kappa$. Let $A$ be the subalgebra of
${\scrp P}(\kappa)$ generated by $[\kappa]^{\leq\rho}$. Then $\cell_{\rm
Hr}(A)=S\cup T\cup U$, where}
\[\eqalignno{S&=\{(\mu,\nu):\omega\leq\mu\leq\nu\leq2^\rho,\nu^\omega=\nu\};\cr
T&=\{(\mu,\mu^\rho):2^\rho<\mu\leq\kappa\};\cr
U&=\{(\mu,\kappa^\rho):2^\rho<\mu,\mu^\rho=\kappa^\rho,\kappa<\mu\}.\cr}\]
Assuming CH and taking $\kappa=\omega_2$ and $\rho=\omega$, we see that
$A$ is the algebra of countable and co-countable subsets of $\omega_2$,
and
\[\cell_{\rm
Hr}(A)=\{(\omega,\omega_1),(\omega_1,\omega_1),(\omega_2,\omega_2)\}.\]
This solves Problem 8(i).

Theorem 10 of Monk, Nyikos [97] says:

\medskip
\par\noindent
{\it Assume GCH, and let ${\scrp A}\subseteq[\kappa^+]^{\kappa^+}$ be such
that the intersection of any two distinct members of ${\scrp A}$ has size
$\leq$ $\kappa$, and $|{\scrp A}|=\kappa^{\kappa^{++}}$. Let $A$ be the
$\kappa^+$-complete subalgebra of ${\scrp P}(\kappa^+)$ generated by
${\scrp A}\cup\{\{\alpha\}:\alpha<\kappa^+\}$. Then}
\[\cell_{\rm
Hr}(A)=\{(\mu,\nu):\omega\leq\mu\leq\nu\leq\kappa^+,
\cf(\nu)>\omega\}\cup\{(\kappa^+,\kappa^{++}),(\kappa^{++},\kappa^{++})\}.\]
Taking $\kappa=\omega$ we get an algebra $A$ such that
\[\cell_{\rm
Hr}(A)=\{(\omega,\omega_1),
(\omega_1,\omega_1),(\omega_1,\omega_2),(\omega_2,\omega_2)\},\]
assuming GCH. This solves part (iii) of Problem 8.

The other two parts of Problem 8 are still open.

\bigskip
\centerline{\bf Problem 9: cellularity for pseudo-tree algebras}

\medskip
\par\noindent
This vague problem is still open.

\bigskip
\centerline{\bf Problem 10: depth of amalgamated free products}

\medskip
\par\noindent
This problem is solved negatively
in Shelah [02], by the following result (see Remark 1.2, 4)):

\medskip
\par\noindent
{\it There is a countable BA $A$ such that for every strong limit cardinal
$\mu$ of cofinality $\aleph_0$ there exist $B,C\geq A$ such that
$\cell(B),\cell(C)\leq\mu$ while $\cell(B\oplus_AC)\geq\mu^+$.}

\bigskip
\centerline{\bf Problem 11: depth of amalgamated free products}

\medskip
\par\noindent
Some consistency results in Shelah [02] partially solve this problem. The
easiest to state is the following part of Observation 1.8:

\medskip
\par\noindent
{\it If $\lambda$ is weakly compact, $A\leq B,C$, $|A|<\lambda$, and
$\depth'(B\oplus_AC)=\lambda^+$, then $\depth'(B)=\lambda^+$ or
$\depth'(C)=\lambda^+$.}

\bigskip
\centerline{\bf Problem 12: depth and ultraproducts}

\medskip
\par\noindent
This problem is solved positively by the following result of Shelah
[$\infty_2$]:

\medskip
\par\noindent
{\it Assume that $\mu$ is a singular cardinal and
$\mu=\mu^\kappa>2^\kappa$. Then there are BAs $B_i$ for $i<\kappa$ such that:

(i) $\depth(B_i)\leq\mu$ for each $i<\kappa$.

(ii) $\mu=\left|\prod_{i<\kappa}\depth(B_i)/D\right|$ for any ultrafilter $D$
on $\kappa$.

(iii) $\depth\left(\prod_{i<\kappa}B_i/D\right)\geq\mu^+$ for any uniform
ultrafilter $D$ on $\kappa$.}

\bigskip
\centerline{\bf Problem 13: tightness and depth}

\medskip
\par\noindent
This problem is partially solved, negatively, by the following result of
Roslanowski, Shelah [00] (Conclusion 7.6):

\medskip
\par\noindent
{\it It is consistent that there is a Boolean algebra $A$
of size $\lambda$
such that there is an ultrafilter of $A$ of tightness $\lambda$, there is
no free sequence of length $\lambda$ in $A$, and
$\t(A)=\lambda\notin\depth_{\rm Hs}(A)$.}

\medskip
\par\noindent
Of course it would be of interest to construct such an example in ZFC.

\bigskip
\centerline{\bf Problem 14: depth and subalgebras}

\medskip
\par\noindent
This problem is partially solved, positively, by the following result of
Roslanowski, Shelah [01], part of Conclusion 18:

\medskip
\par\noindent
{\it It is consistent that there is a cardinal $\kappa$ such that there is
a BA $A$ of size $(2^\kappa)^+$ with $\depth(A)=\kappa$ while
$(\omega,(2^\kappa)^+)\notin\depth_{\rm Sr}(A)$.}

\medskip
\par\noindent
Again it would be of interest to construct such an example in ZFC.

\bigskip
\centerline{\bf Problem 15: depth and subalgebras}

\medskip
\par\noindent
This vague problem is still open.

\bigskip
\centerline{\bf Problem 16: the depth homomorphism relation}

\medskip
\par\noindent
This was solved by Shelah (email message of March 6, 1997). Evidently not
written up yet.

\bigskip
\centerline{\bf Problem 17: the depth homomorphism relation}

\medskip
\par\noindent
This vague problem is still open.

\bigskip
\centerline{\bf Problem 18: topological density and homomorphisms}

\medskip
\par\noindent
This problem was solved by Juhasz and Shelah (email message from Shelah
of March 6,
1997). Evidently not written up yet.

\bigskip
\centerline{\bf Problem 19: $\d_{\rm Hs}$}

\medskip
\par\noindent
This vague problem is still open.

\bigskip
\centerline{\bf Problem 20: $\pi_{{\rm S}+}$}

\medskip
\par\noindent
This problem has an obvious negative solution. For, under
GCH, if $|A|$ is a limit cardinal, then every dense subset of $A$ has size
$|A|$. Suppose now that $|A|=\kappa^+$. Then $\pi_{{\rm S}+}A$ is either
$\kappa$ or $\kappa^+$, and it is clearly attained.

\bigskip
\centerline{\bf Problem 21: $\pi_{\rm Hs}$}

\medskip
\par\noindent
The answer is consistently no according to an email message of Shelah
of March 6, 1997; joint work with Spinas. Evidently this has not been
written up. Of course an example in ZFC
would be interesting.

\medskip
\centerline{\bf Problem 22: length and ultraproducts}

\medskip
\par\noindent
This problem was completely solved by the following result of Shelah
[99]
(Conclusion 15.13 (2)).

\medskip
\par\noindent
{\it  If $D$ is a uniform ultrafilter on $\kappa$,
then for a class of cardinals $\lambda$ such that
$\lambda^\kappa=\lambda$, there is a system $\langle B_i:i<\kappa\rangle$
of Boolean algebras such that $\length(B_i)\leq\lambda$ for each
$i<\kappa$, hence
$\prod_{i<\kappa}\length(B_i)/D)\leq\lambda$, while
$\length(\prod_{i<\kappa}B_i/D)=\lambda^+$.}

\bigskip
\centerline{\bf Problem 23: on $\length_{{\rm h}-}$}

\medskip
\par\noindent
This problem is still open. It may not be clear how this function is
defined; we take its definition to be
\[\length_{{\rm h}-}(A)=\inf\{\sup\{|X|:X\hbox{ is a chain of clopen
subsets of }Y\}:Y\subseteq\ult(A)\}.\]
\vfill
\eject
\centerline{\bf Problem 24: $\irr$ and products}

\medskip
\par\noindent
This was solved positively by Roslanowski, Shelah [00] (Theorem 3.1).

\bigskip
\centerline{\bf Problem 25: $\irr$ and ultraproducts}

\medskip
\par\noindent
A consistent example where this inequality holds was given by Shelah
[$\infty$]. Of course, an example in ZFC would be desirable.

\bigskip
\centerline{\bf Problem 26: $\irr$ and ultraproducts}

\medskip
\par\noindent
A consistent example where this inequality holds was given by Shelah
[99] (part of Conclusion
15.10). Of course, an example in ZFC would be desirable.

\bigskip
\centerline{\bf Problem 27: $\irr$ and $\s$}

\medskip
\par\noindent
Roslanowski and Shelah [00] showed that the answer, consistently, is no,
by the following result, Conclusion 4.6:

\medskip
\par\noindent
{\it It is consistent that there is a BA $B$ such that $\irr(B)=\omega$
and $\s(B\oplus B)=\irr(B\oplus B)=\omega_1$.}

\medskip
\par\noindent
Of course, a counterexample in ZFC would be desirable.

\bigskip
\centerline{\bf Problem 28: $\irr$ and ZFC}

\medskip
\par\noindent
This problem is still open; it is probably difficult.

\bigskip
\centerline{\bf Problem 29: independence and homomorphic images}

\medskip
\par\noindent
This problem is still open.
\bigskip
\centerline{\bf Problem 30: independence and subspaces}

\medskip
\par\noindent
This problem is still open.

\bigskip
\centerline{\bf Problem 31: independence and subspaces}

\medskip
\par\noindent
This problem is still open.

\bigskip
\centerline{\bf Problem 32: independence and cellularity}

\medskip
\par\noindent
Shelah [99] has several results relevant to this problem. A partial
positive solution follows from 6.8 of Shelah [99], which says

\medskip
\par\noindent
{\it If $\kappa$ is weakly inaccessible and
$\langle2^\mu:\mu<\kappa\rangle$ is not eventually constant, then there is
a $\kappa$-cc BA $A$ of size $2^{<\kappa}$ with no independent subset of size
$\kappa^+$.}

\medskip
\par\noindent
This gives, consistently, several examples solving Problem 32
positively. For example, take a model in which $\kappa=\aleph_\alpha$
is weakly inaccessible and $2^{\aleph_\beta}=\aleph_{\alpha+\beta+1}$ for
each regular $\aleph_\beta<\kappa$. Then
$2^{<\kappa}=\aleph_{\alpha+\alpha}$. For Problem 32 one can take
$\rho=\aleph_1$, $\nu=\aleph_3$, and $\lambda=\aleph_{\alpha+3}$.

Of course it would be good to describe completely what happens with
cardinals as in the formulation of Problem 32.

\bigskip
\centerline{\bf Problem 33: independence and cellularity}

\medskip
\par\noindent
Shelah [99] has the following result, which gives a negative solution of
this problem.

\medskip
\par\noindent
(14.24 Conclusion) {\it If $\mu=\mu^{<\mu}<\theta=\theta^{<\theta}$, then
for some $\mu$-complete, $\mu^+$-cc forcing $P$, in $V^P$: If $B$ is a
$\kappa$-cc BA of size $\geq$ $\lambda$, $\mu^{<\kappa}=\mu$, and
$\lambda$ is a regular cardinal in $(\mu,\theta]$, then $\lambda$ is a
free caliber of $B$.}

\medskip
\par\noindent
For Problem 33, we can apply this to a model of GCH, with
$\mu=\aleph_{\omega+1}$ and $\theta=\aleph_{\omega+\omega+1}$. We then get a
model such that if $B$ is an $\aleph_1$-cc BA of size $\geq$
$\aleph_{\omega+\omega+1}$, then $\aleph_{\omega+\omega+1}$
is a free caliber of $B$. This shows
that the case of Problem 33 in which $\mu=\aleph_{\omega+\omega}$,
$\kappa=\aleph_1$, and $\lambda=\aleph_{\omega+\omega+1}$
has a negative solution.

Although this solves Problem 33 as stated, it would be good
to describe completely what happens with
cardinals as in the formulation of Problem 33.

\bigskip
\centerline{\bf Problem 34: independence and cellularity}

\medskip
\par\noindent
The result of Shelah [99] described above for Problem 32 also solves this
problem positively.

\bigskip
\centerline{\bf Problem 35: products and free caliber}

\medskip
\par\noindent
By Shelah [99], 6.11, consistently either answer to this question is
possible.

\bigskip
\centerline{\bf Problem 36: completions and free caliber}

\medskip
\par\noindent
By Shelah [99], 6.11, consistently either answer to this question is
possible.

\bigskip
\centerline{\bf Problem 37: complete Boolean algebras and free caliber}

\medskip
\par\noindent
The comment preceding Problem 37 is wrong; in Monk [83] it is merely
observed that $\hbox{Freecal}(\overline{\hbox{Intalg}L})$
is empty for those
linear orders of size $2^\mu$ with a dense subset of size $\mu$. This easy
observation does imply that it follows from GCH that for every $\mu$ there
is a complete BA of power $2^\mu$ with Freecal empty.

A result in the other direction was obtained by Shelah [99], Claim 8.1:

\medskip
\par\noindent
{\it Assume GCH in the ground model, and let $P$ be the partial order for
adding $\aleph_{\omega_1}$ Cohen reals. Then in the generic extension we
have $2^{\aleph_0}=\aleph_{\omega_1}$, $2^{\aleph_1}=\aleph_{\omega_1+1}$,
and:

There is no complete BA $A$ of size $2^{\aleph_1}$ such that $\hbox{\rm
Freecal}(A)=\emptyset$; in fact, every complete BA of that size has free
caliber $2^{\aleph_1}$.}

\medskip
\par\noindent
Of course this leaves open what happens for other cardinals.

\bigskip
\centerline{\bf Problem 38: specifying the set {\rm Freecal}}

\medskip
\par\noindent
This problem is still open.

\bigskip
\centerline{\bf Problem 39: Interval algebras and finite independence}

\medskip
\par\noindent
This problem is still open.

\vfill
\eject
\centerline{\bf Problem 40: $\pi$ and $\pi\chi$ for complete BAs}

\medskip
\par\noindent
This problem is still open.

\bigskip
\centerline{\bf Problems 41, 42: attainment of tightness}

\medskip
\par\noindent
The following example of J. C. Mart\'\i nez [02] completely solves these
problems negatively, in ZFC;
there was an earlier consistency result of Roslanowski and
Shelah:

\medskip
\par\noindent
{\it Let $\kappa$ be a limit cardinal such that $\lambda\buildrel{\rm
def}\over=\cf\kappa>\omega$. Let
$\langle\kappa_\alpha:\alpha<\lambda\rangle$ be a strictly increasing
sequence of infinite successor cardinals with supremum $\kappa$. For each
$\alpha<\lambda$ let $G_\alpha$ be the ultrafilter on
$\intalg(\kappa_\alpha)$ generated by
$\{[\xi,\kappa_\alpha):\xi<\kappa_\alpha\}$. Let
\[A=\left\{x\in\prod^{\rm
w}_{\alpha<\lambda}\intalg(\kappa_\alpha):
\forall\alpha,\beta<\lambda(x_\alpha\in
G_\alpha\hbox{ \iff }x_\beta\in G_\beta)\right\}\]
Then $A$ is a subalgebra of $\prod^{\rm
w}_{\alpha<\lambda}\intalg(\kappa_\alpha)$. Now let
\[F=\{x\in A:\forall\alpha<\lambda(x_\alpha\in G_\alpha)\}.\]
Then $F$ is an ultrafilter on $A$.

The tightness of $A$ is equal to the tightness of $F$, which is
$\kappa$. On the other hand, $A$ does not have a free sequence of length
$\kappa$.}

\bigskip
\centerline{\bf Problem 43: $\pi\chi$ and attainment of tightness}

\medskip
\par\noindent
This problem is still open.

\bigskip
\centerline{\bf Problem 44: tightness and unions}

\medskip
\par\noindent
This problem has been solved by Roslanowski and Shelah (email message
of March 7, 1997). I do not know where this will appear.

\bigskip
\centerline{\bf Problem 45: tightness and depth}

\medskip
\par\noindent
This problem is solved negatively, consistently, by the following
result of Shelah and Spinas [99]:

\medskip
\par\noindent
{\it Suppose that $0^\sharp$ exists. Let $B$ be a superatomic Boolean
algebra in the constructible universe $L$, and let $\lambda$ be an
uncountable cardinal in $V$. Then in $L$ it is true that
$\t'(B)\geq\lambda^+$ implies that $\depth'(B)\geq\lambda$.}

\medskip
\par\noindent
Here one applies Theorem 6 to get
$R_\omega(\lambda^+,\lambda,\omega)$, and then applies Corollary 5.

Of course, a solution in ZFC would be desirable.

\bigskip
\centerline{\bf Problem 46: ultraproducts and spread}

\medskip
\par\noindent
This problem was solved positively in Shelah [99] by the following
result, part of Conclusion 15.13:

\medskip
\par\noindent
{\it If $D$ is a uniform ultrafilter on $\kappa$, then there is a
class of cardinals $\chi$ with the following properties:

(i) $\chi^\kappa=\chi$.

(ii) There are BAs $B_i$ for $i<\kappa$ such that:

\hskip.2in(a) $\s(B_i)\leq\chi$ for each $i<\kappa$, and hence
$|\prod_{i<\kappa}\s(B_i)/D|\leq\chi$.

\hskip.2in(b) $\s(\prod_{i<\kappa}B_i)=\chi^+$.}

\bigskip
\centerline{\bf Problem 47: ultraproducts and spread}

\medskip
\par\noindent
This was solved positively by the following result of Shelah, Spinas
[$\infty$], part of Corollary 2.4:

\medskip
\par\noindent
{\it There is a model in which there exist cardinals $\kappa,\mu$,
a system $\langle
B_i:i<\kappa\rangle$ of BAs, and an ultrafilter $D$ on $\kappa$ such
that $|\prod_{i<\kappa}\s(B_i)/D|=\mu^{++}$ and
$\s(\prod_{i<\kappa}B_i)\leq\mu^+$.}

\bigskip
\centerline{\bf Problem 48: ultraproducts and character}

\medskip
\par\noindent
This was solved positively by the following result of Shelah, Spinas
[$\infty$], part of Corollary 2.7:

\medskip
\par\noindent
{\it There is a model in which there exist cardinals $\kappa,\mu$,
a system $\langle
B_i:i<\kappa\rangle$ of BAs, and an ultrafilter $D$ on $\kappa$ such
that $|\prod_{i<\kappa}\chi(B_i)/D|=\mu^{++}$ and
$\chi(\prod_{i<\kappa}B_i)\leq\mu^+$.}

\bigskip
\centerline{\bf Problem 49: spread and character}

\medskip
\par\noindent
This problem is still open. It is related to the difficult problem of
general $S$ and $L$ spaces.

\bigskip
\centerline{\bf Problem 50: attainment for $\hl$}

\medskip
\par\noindent
Two main results on this problem are found in Roslanowski, Shelah
[01a]. In section 1 of that paper they show (hypothesis 1.1 and
Theorem 1.4):

\medskip
\par\noindent
{\it Suppose that $\langle\chi_i:i<\cf(\lambda)\rangle$ is a strictly
increasing sequence of infinite cardinals such that $(2^{{\rm
cf}(\lambda)})^+<\chi_0$ and $\lambda=\sup_{i<{\rm
cf}(\lambda)}\chi_i$. Then for any BA $A$ such that $\hl(A)=\lambda$,
all (three) versions of $\hl$ are attained.}

\medskip
\par\noindent
On the other hand, in 3.7 they show that
it is consistent to have a BA with $\hl$ attained in the
right-separated sense but not in the ideal-generation sense.

Recall from [CI] that the three versions of $\hl$ is question are (1)
the version of the definition,
$\sup\{\hbox{L}(X):X\subseteq\ult(A)\}$, (2) the version involving ideal
generation, and (3) the version involving right separated sequences.
Write $\hbox{att}(i)$ to indicate attainment in the sense ($i$) for
$i=1,2,3$.
In [CI] it is shown that $\hbox{att}(2)$ implies $\hbox{att}(3)$, and
the result 3.7 above is that it is consistent to have an example with
$\hbox{att}(3)$ but not $\hbox{att}(2)$. Such an example in ZFC would
be of interest. Examples in ZFC or consistency results for the
remaining implications $\hbox{att}(1)\Rightarrow\hbox{att}(2)$,
$\hbox{att}(1)\Rightarrow\hbox{att}(3)$,
$\hbox{att}(2)\Rightarrow\hbox{att}(1)$,
$\hbox{att}(3)\Rightarrow\hbox{att}(1)$ would be of interest also.

\bigskip
\centerline{\bf Problem 51: $\hl$ and ultraproducts}

\medskip
\par\noindent
This problem was solved in Shelah [99] by the following result, part
of 15.13:

\medskip
\par\noindent
{\it If $D$ is a uniform ultrafilter on $\kappa$ then there is a class
of cardinals $\chi$ with $\chi^\kappa=\chi$ such that there are BAs
$B_i$ for $i<\kappa$ such that:

(i) $\hl(B_i)\leq\chi$, hence $|\prod_{i<\kappa}\hl(B_i)/D|\leq\chi$;

(ii) $\hl(\prod_{i<\kappa}B_i/D)=\chi^+$.}

\bigskip
\centerline{\bf Problem 52: $\hl$ and ultraproducts}

\medskip
\par\noindent
This was solved positively by the following result of Shelah, Spinas
[$\infty$], part of Corollary 2.4:

\medskip
\par\noindent
{\it There is a model in which there exist cardinals $\kappa,\mu$,
a system $\langle
B_i:i<\kappa\rangle$ of BAs, and an ultrafilter $D$ on $\kappa$ such
that $|\prod_{i<\kappa}\hl(B_i)/D|=\mu^{++}$ and
$\hl(\prod_{i<\kappa}B_i/D)\leq\mu^+$.}

\bigskip
\centerline{\bf Problem 53: $\hl$ and $\d$}

\medskip
\par\noindent
This is related to the general $S$ and $L$ space problem, and is
probably difficult.

\bigskip
\centerline{\bf Problem 54: attainment for $\hd$}

\medskip
\par\noindent
In Roslanowski, Shelah [01a] there are two results relevant to this
problem. First, Theorem 1.5 says:

\medskip
\par\noindent
{\it If $2^{{\rm cf}(\lambda)}<\lambda$, then in any BA with $\hd$
equal to $\lambda$, all equivalent versions are attained.}

\medskip
\par\noindent
In section 4, it is shown that it is consistent to have a singular
cardinal $\lambda$ and a BA $A$ with $\hd(A)=\lambda$, with attainment
in the left-separated sense but not in the sense $\kappa_7$ defined in
[CI].

To describe precisely the problems remaining here, let $\hbox{att}(i)$
denote attainment in the following sense, for $i=0,5,7,8$ (keeping to
the notation of [CI]):

\medskip
\par\noindent
$i=0$: the original definition of $\hd$.

\medskip
\par\noindent
$i=5$: the left-separated equivalent.

\medskip
\par\noindent
$i=7$: supremum of $\pi(B)$ for $B$ a homomorphic image of $A$

\medskip
\par\noindent
$i=8$: supremum of $\d(B)$ for $B$ a homomorphic image of $A$

\medskip
\par\noindent
In [CI] the implications $\hbox{att}(8)\Rightarrow
\hbox{att}(0)\Rightarrow
\hbox{att}(7)\Rightarrow
\hbox{att}(5)$ were shown, and the result above is that there is a
consistent example with $\hbox{att}(5)$ but not $\hbox{att}(7)$.

Thus a ZFC example here would be of interest. Consistent examples or
ZFC examples with $\hbox{att}(7)$ but not $\hbox{att}(0)$, and with
$\hbox{att}(0)$ but not $\hbox{att}(8)$ are not known.

\bigskip
\centerline{\bf Problem 55: $\hd$ and ultraproducts}

\medskip
\par\noindent
This problem was solved in Shelah [99] by the following result, part
of 15.13.

\medskip
\par\noindent
{\it Suppose that $D$ is a uniform ultrafilter on $\kappa$. Then there
is a class of cardinals $\chi$ such that $\chi^\kappa=\chi$ and there
are BAs $B_i$ for $i<\kappa$ such that $\hd(B_i)\leq\chi$ for each
$i<\kappa$, hence $\prod_{i<\kappa}\hd(B_i)/D\leq\chi$, while
$\hd(\prod_{i<\kappa}B_i/D)=\chi^+$.}

\vfill
\eject
\centerline{\bf Problem 56: $\hd$ and ultraproducts}

\medskip
\par\noindent
This problem was solved in Shelah, Spinas [$\infty$] by the following
result, part of 2.4:

\medskip
\par\noindent
{\it It is consistent that there exist cardinals $\kappa,\mu$, BAs
$B_i$ for $i<\kappa$, and an ultrafilter $D$ on $\kappa$ such that
$\prod_{i<\kappa}\hd(B_i)/D=\mu^{++}$ while
$\hd(\prod_{i<\kappa}B_i/D)=\mu^+$.}

\bigskip
\centerline{\bf Problems 57, 58: $\hd$, $\s$, and $\chi$}

\medskip
\par\noindent
These two problems are related to the general $S$- and $L$-space
problem, and thus are probably difficult.

\bigskip
\centerline{\bf Problem 59: incomparability and ultraproducts}

\medskip
\par\noindent
This problem is still open.

\bigskip
\centerline{\bf Problem 60: incomparability and ultraproducts}

\medskip
\par\noindent
This problem is solved positively by the following result of Shelah,
Spinas [$\infty$], part of 1.7:

\medskip
\par\noindent
{\it It is consistent that there exist cardinals $\kappa,\mu$ an
ultrafilter $D$ on $\kappa$, and a system
$\langle B_i:i<\kappa\rangle$ of interval algebras
such that $\prod_{i<\kappa}(\inc(B_i)/D)=\mu^{++}$ while
$\inc(\prod_{i<\kappa}B_i/D)\leq\mu^+$.}

\bigskip
\centerline{\bf Problem 61: incomparability and $\chi$}

\medskip
\par\noindent
This problem is still open; it is related to the general $S$- and
$L$-space problem, and is thus probably difficult.

\bigskip
\centerline{\bf Problem 62: $\hcof$ and ultraproducts}

\medskip
\par\noindent
This question is answered positively in Shelah [$\infty$] (see the
comment at the end of the paper).

\bigskip
\centerline{\bf Problem 63: $\hcof$ and $\inc$}

\medskip
\par\noindent
This problem is open. Again it is related to the general $S$- and
$L$-space problems, and so is probably difficult.

\bigskip
\centerline{\bf Problem 64: ideals and automorphisms}

\medskip
\par\noindent
This problem is open.

\bigskip
\centerline{\bf Problem 65: $\irr$ and $\chi$}

\medskip
\par\noindent
This problem is open. It is probably difficult, being related to the
general $S$- and $L$-space problem.

\bigskip
\centerline{\bf Problem 66: $\irr$ and $\inc$}

\medskip
\par\noindent
This problem is open. It is probably difficult.

\bigskip
\centerline{\bf Problem 67: $\irr$ and $\hcof$}

\medskip
\par\noindent
This problem is open. It is probably difficult, being related to the
general $S$- and $L$-space problem.

\bigskip
\centerline{\bf Problem 68: ideals and endomorphisms}

\medskip
\par\noindent
This problem is open.

\bigskip
\centerline{\bf Problem 69: ideals and endomorphisms}

\medskip
\par\noindent
This problem is open.

\bigskip
\centerline{\bf Problem 70: ideals and automorphisms}

\medskip
\par\noindent
This problem is open.

\bigskip
\centerline{\bf Problem 71: endomorphisms for tree algebras}

\medskip
\par\noindent
This problem is open.

\bigskip
\centerline{\bf Problem 72: $\s$ and $\inc$ for superatomic algebras}

\medskip
\par\noindent
This problem is open.

\bigskip
\centerline{\bf Problem 73: $\s$ and $\irr$ for superatomic algebras}

\medskip
\par\noindent
In Roslanowski, Shelah [00] the consistency of a superatomic BA $A$
with $\s(A)<\irr(A)$ is proved (part of Conclusion 5.8).
Of course it would be interesting to
get such an example in ZFC.

\bigskip
\centerline{\bf Problem 74: $\inc$ and cardinality for superatomic
algebras}

\medskip
\par\noindent
This problem is still open.

\bigskip
\centerline{\bf Problem 75: $\irr$ and cardinality for superatomic
algebras}

\medskip
\par\noindent
This problem is still open; it is probably difficult.

\bigskip
\centerline{\bf Problem 76: automorphisms and endomorphisms for
superatomic algebras}

\medskip
\par\noindent
This problem is solved by the following result of Shelah [01],
part of Conclusion 1.8:

\medskip
\par\noindent
{\it Suppose that $\mu$ is a singular strong limit cardinal with
cofinality greater than $\omega$. Also, suppose that $\kappa>\mu$ is
regular and $\kappa\leq2^\mu<2^\kappa$. Then there is a superatomic BA
of size $\kappa$ with $\mu$ atoms, with at most $2^\mu$ automorphisms,
and with $2^\kappa$ endomorphisms.}

\bigskip
\centerline{\bf Problem 77: ideals and subalgebras in superatomic algebras}

\medskip
\par\noindent
This problem is partially solved by the following result in
Roslanowski, Shelah [00] (part of Conclusion 5.8):

\medskip
\par\noindent
{\it It is consistent that for some $\kappa$, there is a superatomic
BA $A$ such that $|\id(A)|=\kappa^+$ and $|\sub(A)|=2^{\kappa^+}$.}

\medskip
\par\noindent
Of course, an example in ZFC would be of interest.

\bigskip
\centerline{\bf Problem 78: $\inc$ and $\irr$ in superatomic algebras}

\medskip
\par\noindent
This problem is partially solved by the following result in
Roslanowski, Shelah [00] (part of Conclusion 5.8):

\medskip
\par\noindent
{\it It is consistent that for some $\kappa$, there is a superatomic
BA $A$ such that $\inc(A)=\kappa$ and $\irr(A)=\kappa^+$.}

\medskip
\par\noindent
Of course, an example in ZFC would be of interest.

\bigskip
\centerline{\bf Problem 79: $\inc$ and $\irr$ in superatomic algebras}

\medskip
\par\noindent
This problem is partially solved by the following result in
Roslanowski, Shelah [00] (conclusion 6.5):

\medskip
\par\noindent
{\it It is consistent that for some $\kappa$, there is a superatomic
BA $A$ such that $\inc(A)=\kappa^+$ and $\irr(A)=\kappa$.}

\medskip
\par\noindent
Of course, an example in ZFC would be of interest, but this is
probably difficult.

\bigskip
\centerline{\bf Problem 80: automorphisms in superatomic algebras}

\medskip
\par\noindent
This problem is solved by the following result of Shelah [01], part of
Theorem 2.2:

\medskip
\par\noindent
{\it If $\mu$ is a singular strong limit cardinal of cofinality
$\omega$, then there is a superatomic BA $A$ such that
$|\aut(A)|\leq2^\mu<|A|$.}

\bigskip
\centerline{\bf Problem 81: automorphisms and $\t$}

\medskip
\par\noindent
This problem is solved by the following result of Roslanowski, Shelah
[00] (Conclusion 6.10):

\medskip
\par\noindent
{\it It is consistent that for some cardinal $\kappa$ there is a
superatomic BA $A$ such that $|\aut(A)|=\kappa$ and $\t(A)=\kappa^+$.}

\medskip
\par\noindent
However, it would be of interest to get such an example in ZFC.

\bigskip
\centerline{\bf Problem 82: $\s$ and $\hl$ for atomic algebras}

\medskip
\par\noindent
This problem is still open, and is related to the difficult general
$S$- and $L$-space problem.

\bigskip
\centerline{\bf Problem 83: $\s$ and $\hd$ for atomic algebras}

\medskip
\par\noindent
This problem is still open, and is related to the difficult general
$S$- and $L$-space problem.

\bigskip
\centerline{\bf Problems 84--91: semigroup algebras}

\medskip
\par\noindent
There does not seem to have been any work on these problems. Problems
85, 87, and 91 appear to be difficult.

\bigskip
\centerline{\bf Problem 92: automorphisms and $\pi\chi$ for semigroup
algebras}

\medskip
\par\noindent
This problem has an obvious positive answer: take a rigid
interval algebra. One can have such an algebra with $\pi\chi$
arbitrarily large. See, for example, Monk [96].

\bigskip
\centerline{\bf Problem 93: automorphisms and $\ind$ in semigroup
algebras}

\medskip
\par\noindent
This problem also has an obvious positive answer with a
rigid interval algebra. It is still open to find a rigid semigroup
algebra with $\ind$ arbitrarily large.

\bigskip
\centerline{\bf Problems 94--97: semigroup algebras and algebras}
\centerline{\bf of Kunen and of Baumgartner, Komjath}

\medskip
\par\noindent
No work has been done on these problems.

\bigskip
\centerline{\sectbf Bibliography}

\medskip
\par\noindent
Juh\'asz, I.; Shelah, S. [98] {\it On the cardinality and weight spectra
of compact spaces, II.} Fund. Math. 155 (1998), 91--94. Publication 612 of
Shelah.

\medskip
\par\noindent
Koppelberg, S. [89] {\bf General theory of Boolean algebras.} Vol. 1 of
{\bf Handbook of Boolean algebras.} North-Holland. 312pp.

\medskip
\par\noindent
Magidor, M.; Shelah, S. [98] {\it Length of Boolean algebras and
ultraproducts.} Math. Japon. 48 (1998), 301--307. Publiction 433 of Shelah.

\medskip
\par\noindent
Martinez, J. C. [02] {\it Attainment of tightness in Boolean spaces.}
Math. Log. Quart. 48 (2002), no. 4, 555--558.

\medskip
\par\noindent
Monk, J. D. [CI] {\bf Cardinal 
invariants on Boolean algebras.} (1996) Birkh\"auser
Verlag, 298pp.

\medskip
\par\noindent
Monk, J. D. [96] {\it Minimum-sized infinite partitions of Boolean algebras.}
Math. Logic Quarterly 42 (1996), 537--550.

\medskip
\par\noindent
Monk, J. D.; Nyikos, P. [97] {\it On cellularity in homomorphic images of
Boolean algebras.} Topol. Proc. 22, summer 1997, 341--362.

\medskip
\par\noindent
Roslanowski, A.; Shelah, S. [00] {\it More on cardinal invariants of
Boolean algebras.} Ann. Pure Appl. Logic 103 (2000), 1--37. Publication
599 of Shelah.

\medskip
\par\noindent
Roslanowski, A.; Shelah, S. [01] {\it Historic forcing for {\rm Depth}.}
Colloq. Math. 89 (2001), 99--115. Publication 733 of Shelah.

\medskip
\par\noindent
Roslanowski, A.; Shelah, S. [01a] {\it Forcing for $\hl$ and $\hd$}
Colloq. Math. 88 (2001), 273-310. Publication 651 of Shelah.

\medskip
\par\noindent
Shelah, S. [91] {\it Strong negative partition relations below the
continuum.} Acta Math. Hung. 58, no. 1-2 (1991), 95--100. Publication
no. 327.

\medskip
\par\noindent
Shelah, S. [94] {\it $\aleph_{\omega+1}$ has a J\'onsson algebra.} Chapter
II of {\bf Cardinal Arithmetic.} Oxford Univ. Press, 34--116. Publication
no. 355.

\medskip
\par\noindent
Shelah, S. [94$'$] {\it J\'onsson algebras in inaccessible cardinals.}
Chapter III of
{\bf Cardinal Arithmetic.} Oxford Univ. Press, 34--116. Publication
no. 365.

\medskip
\par\noindent
Shelah, S. [94$''$] {\it Colorings.} Appendix I of
{\bf Cardinal Arithmetic.} Oxford Univ. Press, 34--116. Publication
no. 282a.

\medskip
\par\noindent
Shelah, S. [96] {\it Was Sierpi\'nski right III? Can continuum-c.c. times
c.c.c be continuum-c.c.?} Annals of Pure and Applied Logic 78 (1996),
259--269. Publication no. 481.

\medskip
\par\noindent
Shelah, S. [97] {\it Colouring and non-productivity of $\aleph_2$-c.c.}
Annals of Pure Appl. Logic 84, issue 2 (1997), 153-174. Publication no.
572.

\medskip
\par\noindent
Shelah, S. [99] {\it Special subsets of $^{{\rm cf}(\mu)}\mu$, Boolean
algebras and Maharam measure algebras.} Topol. Appl. 99 (1999), 135--235.
Publication 620.

\medskip
\par\noindent
Shelah, S. [00] {\it Cellularity of free products of Boolean algebras.}
Fundamenta Math. 166 (2000), 153-208. Publication no. 575.

\medskip
\par\noindent
Shelah, S. [00a] {\it Applications of PCF theory.} J. Symb. Logic 65
(2000), 1624--1674. Publication no. 589.

\medskip
\par\noindent
Shelah, S. [01] {\it Constructing Boolean algebras for cardinal
invariants.} Alg. Univ. 45 (2001), 353--373. Publication no. 641.

\medskip
\par\noindent
Shelah, S. [02] {\it More constructions for Boolean algebras.}
Arch. Math. Logic 41 (2002), no. 5, 401-441. Publication no. 652.

\medskip
\par\noindent
Shelah, S. [$\infty$] {\it On ultraproducts of Boolean algebras and
$\irr$.}
Publication no. 703.

\medskip
\par\noindent
Shelah, S. [$\infty_2$] {\it On Monk's problems: the depth of ultraproducts.}
Publication no. 645.

\medskip
\par\noindent
Shelah, S.; Spinas, O. [99] {\it On tightness and depth in superatomic
Boolean algebras.} Proc. Amer. Math. Soc. 127 (1999),
3475--3480. Publication 663 of Shelah.

\medskip
\par\noindent
Shelah, S.; Spinas, O. [$\infty$] {\it On incomparability and related
cardinal functions on ultraproducts of Boolean algebras.} Publication
677 of Shelah.

\bye